\documentclass[12pt,twoside,doublespace]{article}

\usepackage{amssymb}
\usepackage{graphicx}

\def\smskip{\par\vskip 5 pt}
\def\QED{\hfill $\Box$\smskip}
\newtheorem{theorem}{Theorem}
\newtheorem{lemma}{Lemma}
\newtheorem{proposition}{Proposition}

\begin{document}

\begin{center}

\vspace{35pt}

{\Large \bf Selective Bi-coordinate Method}

\vspace{5pt}

{\Large \bf for Non-Stationary and Non-Smooth }

\vspace{5pt}

{\Large \bf Resource Allocation Type Problems}

\vspace{35pt}

{\sc I.V.~Konnov\footnote{\normalsize E-mail: konn-igor@ya.ru}}

\vspace{35pt}

{\em  Department of System Analysis
and Information Technologies, \\ Kazan Federal University, ul.
Kremlevskaya, 18, Kazan 420008, Russia.}

\end{center}

\vspace{35pt}

\begin{abstract}
We propose a method of bi-coordinate variations for non-stationary and non-smooth
optimization problems, which involve a single linear equality and box
constraints. Here only approximation sequences are known instead of exact values
of the cost function and parameters of the feasible set.
It consists in making descent steps with respect to only two selected
coordinates satisfying some special threshold rule.
The method is simpler essentially than the usual gradient or dual type ones
and differs from the previous known bi-coordinate ones suggested for the usual
stationary and smooth problems. We establish its convergence under rather mild assumptions.
Computational tests also reveal certain preferences of the proposed method over the
known ones.

{\bf Key words:} Optimization problems, non-stationary, non-smooth functions, linear equality
constraint, box constraints, bi-coordinate variations, threshold control.
\end{abstract}

{\bf MSC codes:}{ 90C30, 90C25, 90C06, 91B32, 68W15}

\newpage


\section{Introduction}\label{s1}

The custom finite-dimensional optimization problem consists in finding
the minimal value of some goal function $f : \mathbb{R}^{n} \to \mathbb{R}$ on a feasible set $D$
such that $D  \subseteq \mathbb{R}^{n}$. For brevity, we
write this problem as
\begin{equation} \label{eq:1.1}
 \min \limits _{x \in D} \to f(x),
\end{equation}
its solution set is denoted by $D^{*}$ and the optimal value of the
function by $f^{*}$, i.e. $f^{*} = \inf\limits_{x \in D} f(x)$.
Many problems of optimal allocation of some resource within a composite system
containing $n$ elements  can be reduced to the above format
where
\begin{equation} \label{eq:1.2}
D=\left\{ x\in X \ \vrule \ \langle a,x\rangle=\beta
\right\} \ \mbox{and} \ X=[\alpha' _{1},\alpha'' _{1}]\times \ldots \times [\alpha' _{n},\alpha'' _{n}],
\end{equation}
$\beta$ is a fixed number, $a=(a_{1}, \ldots, a_{n})^{\top}$ is a fixed
vector whose coordinates are non-zero, 
and $\langle c,d\rangle$ denotes the usual scalar
product of vectors $c$ and $d$; see e.g.
\cite{BT89}--\cite{Pat08} and references therein.
Then, solution of problem (\ref{eq:1.1})--(\ref{eq:1.2}) yields
a feasible resource allocation that minimizes the total system dis-utility.
Such problems arise in various fields and are investigated rather well and many rather
efficient algorithms have been proposed; see e.g.
\cite{BT89,Pat08} and references therein.

However, the recent development of communication and information processing
technologies reveal special features of resource allocation problems arising
in these fields; see e.g. \cite{CW03,SWB06} and references therein.
Namely,  they also reduce to the form (\ref{eq:1.1})--(\ref{eq:1.2}), but have
very large dimensionality, inexact and/or non-stationary parameters reflecting
variability of users' behavior, and scattered necessary information.
Hence, we are forced to develop methods whose iteration computation
expenses and accuracy requirements are rather low and
do not utilize matrix transformations at each iteration as
the Newton or interior point type ones. This means that
even simple coordinate-wise descent methods
may appear very useful here.

Besides, the same
optimization formulation is paid now a significant attention  due to its
various big data applications; see e.g. \cite{Bur98,CBS14}
and the references therein. In fact, similar optimization problems arise in
machine learning, signal, speech and image recognition and processing,
and so on. These problems possess almost the same features; i.e.,
huge dimensionality, inexact,  incomplete, and/or non-stationary data,
which can be scattered within different computer networks.
 Moreover, they are often contain non-smooth regularization
or penalty terms and rather simple constraints. As a result,
even calculation of all the components of the
gradient may be too hard. This fact creates certain difficulties for application of
 custom second and even first order optimization methods.

For this reason, we are interested in
developing special low cost iterative methods, which are applicable for
problems of form (\ref{eq:1.1})--(\ref{eq:1.2}) and keep the convergence
properties of the usual ones. In particular, their  computational
expenses per iteration should be reduced essentially.

In this paper, we intend to develop a new bi-coordinate descent method
for these problems. We recall that the first bi-coordinate method
for problems with one linear constraint and only lower (non-negativity)
bounds for variables was proposed and substantiated in
\cite{Kor80}. Further, these methods became very popular due to
their rather good performance for data mining applications; see e.g. \cite{Pla99}--\cite{GI06};
the detailed description of the recent versions is given e.g.
in \cite{LPR09,Bec14}.

However, most of these methods are based on either computation of
certain marginal indices or utilization of some general error bound
and Lipschitz constants for the gradient, so that finding a descent
direction in these methods will require calculation of all the partial derivatives at
each iteration, i.e., their iteration cost is almost the same as in the usual projection or
conditional gradient methods. The other methods exploit
the random coordinate choice idea, which reduces computational expenses per iteration,
but may however lead to rather slow convergence.

Rather recently, a so-called selective bi-coordinate method with special
threshold control and tolerances was proposed in \cite{Kon16a} for
problem (\ref{eq:1.1})--(\ref{eq:1.2}) with $\alpha' _{i}=0$ and $\alpha'' _{i}=\infty$
for all $i =1, \ldots, n$, besides, the vector $a$ was chosen to be the vector of units,
that is, it was destined for the case of the usual simplex constraints.
Its bi-coordinate descent is based on satisfying
some threshold value and does not require calculation of all the
partial derivatives in general. Besides, its threshold
control strategy seems more flexible in comparison with the previous rules.
In \cite{Kon16}, its complexity estimate
 $O(1/\alpha)$, which gives the the total number of
iterations for attaining the accuracy $\alpha$, was established
for the case where the goal function is convex
and its partial bi-coordinate gradients are Lipschitz continuous.
It should be noted that this method can be treated  as a
self-adjustment process for attaining an equilibrium state of
a closed economic system; see \cite{Kon16a,Kon15e}. However, this method
can not be applied directly to a general problem of form (\ref{eq:1.1})--(\ref{eq:1.2})
with both upper and lower bounds having different signs, which
somewhat restricts its field of applications.

The main goal of this paper is to develop a new selective bi-coordinate method,
which follows the approach from \cite{Kon16a}, but becomes suitable for general
non-stationary and non-smooth optimization problems of form (\ref{eq:1.1})--(\ref{eq:1.2}).
This means that only approximation sequences are known instead of exact values
of the cost function and parameters of the feasible set, besides, the
limit goal function $f$ can be non-smooth. Clearly, these properties enlarges
its areas of significant applications essentially. We establish its convergence
and report some results of computational experiments
with the new method and compare them with some related ones.


\section{Basic preliminaries and assumptions}\label{s2}

We will use the following first set of basic
 assumptions for problem (\ref{eq:1.1})--(\ref{eq:1.2}).

{\bf (A1)}  {\em The feasible set  $D$ is nonempty, the set  $X$ is bounded,
$a_{i}>0$ for all $i \in I=\{1, \ldots, n\}$,
the function $f : \mathbb{R}^{n} \to \mathbb{R}$
is locally Lipschitz on $X$, i.e. it is Lipschitz continuous in a neighborhood of
any point $x \in X$.}

Then problem (\ref{eq:1.1})--(\ref{eq:1.2}) has a
solution and $f^{*}>- \infty$.
Here we notice that the positivity of $a_{i}$ does not restrict the generality.
In fact, if all $a_{i}$ are negative, we can obtain the previous case by simple replacing
$\beta$ with $-\beta$. Next, we can consider a somewhat more general case
where $a$ has arbitrary non-zero coordinates.
However, since the signs of lower and upper bounds are also arbitrary, we can
introduce the new variables $y_{i}={\rm sign} (a_{i})x_{i}$ for all $i \in I$
and insert the new bounds  $\tilde \alpha' _{i}=-\alpha'' _{i}$, $\tilde \alpha'' _{i}=-\alpha' _{i}$ if
$a_{i}<0$ together with the previous ones $\tilde \alpha' _{i}=\alpha' _{i}$, $\tilde \alpha'' _{i}=\alpha'' _{i}$ if
$a_{i}>0$. In such a way,  we again obtain the problem of form (\ref{eq:1.1})--(\ref{eq:1.2})
satisfying the above assumptions; see also
\cite{Bec14}. In the other words, we can always obtain the same sign for all the entries of $a$ by
proper changes of lower and upper bounds of variables.

We now recall some concepts and properties from Non-smooth Analysis; see \cite{Cla83}
for more details. Since $f$ is Lipschitz continuous in a neighborhood of
$x  \in X$, we can define its generalized gradient set at $x$:
$$
\partial^{\uparrow} f(x)=\{ g \in \mathbb{R}^{n} \ | \
\langle g , p \rangle \le f^{\uparrow}(x;p)\},
$$
which must be non-empty, convex and closed. Here $f^{\uparrow}(x;p)$ denotes the upper Clarke-Rockafellar
derivative:
$$
  f ^{\uparrow}(x;p) = \limsup _{y \to x, \alpha \searrow 0}
 ((f (y+ \alpha p) - f (y))/ \alpha ).
$$
It follows that
$$
  f ^{\uparrow}(x,p) = \sup _{g \in \partial^{\uparrow} f(x)}
            \langle g, p \rangle.
$$
At the same time,  the
function  $f$ has the gradient $\nabla f(x)$ a.e. in $X$,
furthermore, it holds that
\begin{equation} \label{eq:2.1}
  \partial^{\uparrow} f(x)
={\rm conv} \left\{\lim \limits_{y\rightarrow x}\nabla f(y)\ | \
y\in D_{f}, \ y\notin S \right\},
\end{equation}
where $D_{f}$ denotes  the set of points where $f$ is
differentiable, and $S$ denotes  an arbitrary subset of measure
zero. If $f$ is convex, then  $\partial^{\uparrow} f(x)$ coincides
with the subdifferential $\partial f(x)$ in the sense of Convex Analysis,
i.e.,
$$
 \partial f(x) = \{  g \in \mathbb{R}^{n}  \ | \ f(y)-
  f(x) \geq  \langle g, y-x \rangle   \quad \forall y \in \mathbb{R}^{n} \}.
$$
In this case, we have
$$
  f'(x;p) = \lim_{\alpha \to 0}
      ((f (x+ \alpha p) - f (x))/ \alpha )
         = \sup _{g \in \partial f(x)}
            \langle g, p \rangle
$$
and the upper
derivative coincides with the usual direction derivative:
\begin{equation} \label{eq:2.2}
  f ^{\uparrow}(x;p) = f'(x;p).
\end{equation}
Also, if $f$ is differentiable at $x$, (\ref{eq:2.2}) obviously
holds and we have
$$
  f'(x;p) = \langle  \nabla f(x), p \rangle  \
   {\rm and } \  \partial ^{\uparrow} f(x) = \{ \nabla f(x) \};
$$
cf. (\ref{eq:2.1}).

We recall that a function $\varphi : \mathbb{R}^{n}
\to \mathbb{R}$ is called

(a) {\em pseudo-convex} on a set $X$, if
for each pair of points $x, y \in X$, we have
$$
    \varphi' (x; y - x ) \geq 0 \
      \Longrightarrow \
    \varphi (y) \geq \varphi (x);
$$

(b) {\em semi-convex} (or {\em upper pseudo-convex}) if
for each pair of points $x, y \in X$, we have
$$
    \varphi^{\uparrow}(x; y - x ) \geq 0 \
      \Longrightarrow \
    \varphi (y) \geq \varphi (x);
$$
see \cite{Mif77} and also \cite{Kon13}. In case (\ref{eq:2.2}), these concepts coincide, but
in general (b) implies (a). Besides, the class of convex functions is strictly contained
in that of pseudo-convex functions.
We now recall the known optimality condition; see e.g. \cite{Cla83,Mif77} and  \cite{Kon13}.


\begin{proposition} \label{pro:2.1}

(a)  Each solution of problem (\ref{eq:1.1})--(\ref{eq:1.2}) is a solution of the
variational inequality (VI for short): Find a point $x^{*} \in D$
such that
\begin{equation} \label{eq:2.3}
\exists g^{*} \in \partial^{\uparrow} f (x^{*}),  \quad  \langle g^{*},x-x^{*} \rangle \geq 0 \quad \forall x \in D.
\end{equation}

(b)  If $f$ is semi-convex, then  each solution of VI
(\ref{eq:2.3}) solves  problem (\ref{eq:1.1})--(\ref{eq:1.2}).
\end{proposition}

Solutions of VI (\ref{eq:2.3}) are called stationary points
of (\ref{eq:1.1}). It will be suitable to
specialize optimality conditions for the constraints in
(\ref{eq:1.2}).


\begin{proposition} \label{pro:2.2} A point $x^{*}$ is a solution of
VI (\ref{eq:2.3}), (\ref{eq:1.2}) if and only if it satisfies each
of the following equivalent conditions:
\begin{eqnarray}
&& x^{*} \in D, \exists g^{*} \in \partial^{\uparrow} f (x^{*}), \ \exists \lambda, \
     \langle g^{*}-\lambda a,x-x^{*} \rangle \geq 0 \quad \forall x \in X; \label{eq:2.4}\\
&& x^{*} \in D, \exists g^{*} \in \partial^{\uparrow} f (x^{*}), \ \exists \lambda, \ (1/a_{i}) g_{i}^{*} \left\{ {
\begin{array}{ll}
\displaystyle
\geq \lambda \quad & \mbox{if} \ x^{*}_{i}=\alpha' _{i}, \\
=\lambda \quad & \mbox{if} \ x^{*}_{i}\in (\alpha'_{i},\alpha ''_{i}), \\
\leq \lambda \quad & \mbox{if} \ x^{*}_{i}=\alpha'' _{i},
\end{array}
} \right. \quad \mbox{for} \quad i \in I; \label{eq:2.5}\\
&&
\begin{array}{ll}
x^{*} \in D, \ \exists g^{*} \in \partial^{\uparrow} f (x^{*}), & \ \forall i,j \in I, i \neq j,  \\
 & (1/a_{i}) g_{i}^{*}>(1/a_{j}) g_{j}^{*} \ \Longrightarrow \
x^{*}_{i}=\alpha' _{i} \ \mbox{or} \ x^{*}_{j}=\alpha'' _{j};
\end{array} \label{eq:2.6}\\
&&
\begin{array}{ll}
x^{*} \in D, \ \exists g^{*} \in \partial^{\uparrow} f (x^{*}), & \ \forall i,j \in I, i \neq j, \\
& x^{*}_{i} \in (\alpha '_{i},\alpha ''_{i}],  x^{*}_{j} \in [\alpha
'_{j},\alpha ''_{j}) \ \Longrightarrow
\ (1/a_{i}) g_{i}^{*}\leq (1/a_{j}) g_{j}^{*}.
\end{array} \label{eq:2.7}
\end{eqnarray}
\end{proposition}
{\bf Proof.}
In fact, equivalence of (\ref{eq:2.3}) and (\ref{eq:2.4}) follows
from the usual optimality conditions for VIs; see e.g.
\cite[Theorem 12.3]{Kon13}. The equivalence of (\ref{eq:2.4}) and (\ref{eq:2.5}) is
obvious; see e.g. \cite[Proposition 7.2]{Kon07}.

For brevity, set $h_{i}=(1/a_{i}) g_{i}^{*}$ and $h_{j}=(1/a_{j}) g_{j}^{*}$.
Let now a  point $x^{*} \in D$ satisfy
(\ref{eq:2.5}). If there exist $i,j \in I$, $i \neq j$ such that
$h_{i}>h_{j}$, $x^{*}_{i}>\alpha' _{i}$,  and $x^{*}_{j}<\alpha'' _{j}$, then
$h_{i}\leq\lambda$ and $h_{j}\geq\lambda$, which is a contradiction. Hence,
(\ref{eq:2.5}) implies (\ref{eq:2.6}). Clearly, (\ref{eq:2.6})
implies (\ref{eq:2.7}). Let now a  point $x^{*} \in D$ satisfy
(\ref{eq:2.7}). Define the index sets: $I_{-}=\{ i \in I \ | \ x^{*}_{i}=\alpha
'_{i} \} $, $I_{0}=\{ i \in I \ | \ x^{*}_{i}\in (\alpha '_{i},
\alpha ''_{i}) \} $, $I_{+}=\{ i \in I \ | \ x^{*}_{i}=\alpha
''_{i} \} $.

If $I_{0}\neq\varnothing$,  set $\lambda =h_{s}$ for some $s \in I_{0}$.
Then $\lambda=h_{i}$ for any $i \in I_{0}$,
$\lambda \leq h_{i}$ for any $i \in I_{-}$, and
$\lambda\geq h_{i}$ for any $i \in I_{+}$
due to (\ref{eq:2.7}), hence (\ref{eq:2.5}) holds.

Let now $I_{0}=\varnothing$. Set
$\tau_{1}=\max_{ i \in I_{+} } h_{i} $ and
$\tau_{2}=\min_{ i \in I_{-} } h_{i} $, then
(\ref{eq:2.7}) gives $\tau_{1}\leq \tau_{2}$. Take any number
$ \lambda  \in [\tau_{1},\tau_{2}]$, then, by definition,
$\lambda \leq h_{i}$ for any $i \in I_{-}$, and
$\lambda\geq h_{i}$ for any $i \in I_{+}$, which also yields (\ref{eq:2.5}). \QED

We intend to consider the case of the non-stationary optimization problem, where
only sequences of approximations are known instead of the exact
values. This means that we have some sequence of problems of the form:
\begin{equation} \label{eq:2.8}
 \min \limits _{x \in D_{l}} \to f_{l}(x),
\end{equation}
where
\begin{equation} \label{eq:2.9}
D_{l}=\left\{ x\in X_{l} \ \vrule \ \langle a^{l},x\rangle=\beta_{l}
\right\} \ \mbox{and} \ X_{l}=[\alpha' _{1l},\alpha'' _{1l}]\times \ldots \times [\alpha' _{nl},\alpha'' _{nl}],
\end{equation}
$\beta_{l}$ is a fixed number, $a^{l}=(a^{l}_{1}, \ldots, a^{l}_{n})^{\top}$ is a fixed
vector, for $l=0,1,2,\ldots$ The basic approximation assumptions are the following.

{\bf (A2)} {\em  For each $l=0,1,2,\ldots$, the set $D_{l}$ is nonempty, $a^{l}_{i}>0$ and $-\infty < \alpha' _{il} < \alpha'' _{il} <+\infty$ for
$i=1,\ldots,n$,}
\begin{eqnarray*}
   && \lim \limits _{l \to \infty} \alpha' _{il}=\alpha' _{i}, \lim \limits _{l \to \infty} \alpha'' _{il}=\alpha'' _{il}\ \mbox{for} \
i=1,\ldots,n; \\
 && \lim \limits _{l \to \infty} a^{l}=a, \lim \limits _{l \to \infty} \beta_{l}=\beta. \\
\end{eqnarray*}

{\bf (A3)} {\em Each function $f_{l}: X_{l} \rightarrow
\mathbb{R}$ is smooth, the relations $\{ y^{l} \} \to \bar y$  and $y^{l} \in D_{l}$
imply $\{f'_{l}(y^{l}) \} \to \bar g \in \partial^{\uparrow}f(\bar y)$.}

Assumption {\bf (A2)} means that the sequence of the sets $\{ D_{l}\}  $
converges to the limit feasible set $D$, whereas {\bf (A3)} determines
some convergence property of the sequence of the differentiable functions $f_{l}$ to the
non-differentiable function $f$. These assumptions do not seem too restrictive
because they do not include evaluation and concordance of deviations.
In fact, {\bf (A3)} may be invoked by several circumstances.
Firstly, the limit function $f$ of the initial problem may be smooth, and we
replace it with more suitable approximations (say, if $f$ is only convex,
we can take $f_{l}$ strongly convex) or $f_{l}$ remains smooth despite the
inexact calculation of coefficients of $f$. Next, if $f$ is non-smooth,
we can replace it with its smooth approximations. This technique is
well known; see e.g. \cite{ENW95}--\cite{BT12}.  Since $f$ is
locally Lipschitz, it is easy to find such an approximation satisfying
{\bf (A3)}; see \cite{ENW95,CR06}.
There are simple examples for most popular non-smooth functions.
For instance, we can replace $|\tau|$ with $\mu_{1} (\tau,\varepsilon)=\sqrt{\tau^{2}+\varepsilon}$ or
$$
\mu_{2} (\tau,\varepsilon)=\left\{ {
\begin{array}{ll}
\displaystyle
\tau^{2}/2 \quad & \mbox{if} \ |\tau| \leq \varepsilon, \\
\varepsilon \tau -\varepsilon^{2}/2 \quad & \mbox{if} \ |\tau| > \varepsilon;
\end{array}
} \right.
$$
where $\varepsilon>0$ is an approximation parameter.
Nevertheless, we can take into account all the opportunities mentioned above
for approximating the goal function $f$ in order
to enhance the solution method performance.


\section{Some examples of applications}\label{s3}

We intend now to give some examples of applied problems
which reduce to an optimization problem
of form (\ref{eq:1.1})--(\ref{eq:1.2}) and satisfy the above assumptions.


\subsection{Data classification problems}

One of the most popular approaches to data classification is
support vector machine techniques; see e.g. \cite{Bur98,Aga15}.
The simplest linear support vector machine problem for binary data classification
consists in finding a hyperplane separating two collections of known
points $b^{i} \in \mathbb{R}^{m}$, $i=1, \ldots,l$ attributed to some observations
with different labels $ \gamma_{i} \in \{-1,+1\}$, $i=1, \ldots,l$,
where $m$ is the number of features. That is, the distance between the hyperplane and each
collection should be as long as possible. This separation of the
feature space enables us to classify new data points.
However, this requirement appears too strong for real problems
where the so-called soft margin approach, which minimizes the penalties for
mis-classification, is utilized. This problem can be formulated
as the optimization problem
$$
 \min \limits _{w \in \mathbb{R}^{n}} \to (1/p)\|w\|^{p}_{p}
 + C \sum \limits^{l} _{i=1} L( \langle w, b^{i} \rangle-\beta;\gamma_{i})^{q},
$$
where $L$ is a loss function and $C > 0$ is a penalty parameter.
The custom choice is $L(z; y) = \max\{0; 1-yz\}$ whereas
$p$ and $q$ are  either 1 or 2. The more usual 2-norm provides
so useful properties as smoothness of cost functions and uniqueness of solution, but
the 1-norm approach (see \cite{CV95,BO98})
is also very popular since it yields sparsity, i.e. only few
solution components appear non-zero. Due to very large dimensionality
of the feature space, this property is valuable. If
we take  $p=q=1$, we can rewrite this problem as
$$
 \min \limits _{w, \xi} \to \sum \limits^{m} _{j=1} |w_{j}|
 + C \sum \limits^{l} _{i=1} \xi_{i},
$$
subject to
$$
 1+\gamma_{i} (\beta-\langle w, b^{i} \rangle) \leq \xi_{i}, \ \xi_{i} \geq 0, \ i=1, \ldots,l;
$$
or in the equivalent linear programming format:
\begin{equation} \label{eq:3.1}
\min \limits _{u,v, \xi} \to \sum \limits^{m} _{j=1} (u_{j}+v_{j})
 + C \sum \limits^{l} _{i=1} \xi_{i},
\end{equation}
subject to
\begin{eqnarray*}
 && 1+\gamma_{i}\left\{\beta- \sum \limits^{m} _{j=1}(u_{j}-v_{j})b^{i}_{j}\right\}  \leq \xi_{i}, \ i=1, \ldots,l;\\
  && u_{j} \geq 0, v_{j} \geq 0, \ j=1, \ldots,m; \ \xi_{i} \geq 0, \ i=1, \ldots,l;
\end{eqnarray*}
where $w_{j}=u_{j}-v_{j}$, $u_{j} \geq 0$, $v_{j} \geq 0$, and $|w_{j}|=u_{j}+v_{j}$.
We can write now its dual formulation:
\begin{equation} \label{eq:3.2}
\max \limits _{y} \to \sum \limits^{l} _{i=1} y_{i},
\end{equation}
subject to
\begin{eqnarray*}
&& -1 \leq \sum \limits^{l} _{i=1}a_{ij} y_{i} \leq 1, \ i=1, \ldots,l;\\
 && \sum \limits^{l} _{i=1}\gamma_{i} y_{i}=0, \\
  && \ y_{i} \geq 0, \ i=1, \ldots,l;
\end{eqnarray*}
where $a_{ij}=\gamma_{i} b^{i}_{j}$. Utilization just (\ref{eq:3.2}) instead of (\ref{eq:3.1})
is suitable if $l \ll m$, moreover, (\ref{eq:3.2}) allows one to insert new data observations
by simple adding new zero variables without losing the feasibility of the current point.
It seems also worthwhile to replace the first series of double inequalities
with the corresponding penalty:
$$
\min \limits _{y} \to   (\tau/p) \sum \limits^{m} _{j=1}
\left\{ \left( \sum \limits^{l} _{i=1}a_{ij} y_{i} - 1  \right)_{+}^{p}
 + \left(  - \sum \limits^{l} _{i=1}a_{ij} y_{i} - 1  \right)_{+}^{p}\right\}-\sum \limits^{l} _{i=1} y_{i},
$$
subject to
$$
\sum \limits^{l} _{i=1}\gamma_{i} y_{i}=0,  \ y_{i} \geq 0, \ i=1, \ldots,l;
$$
with $\tau>0$, $p$ is  either 1 or 2, $(a)_{+}=\max\{a,0\}$. Clearly, this problem falls into format
(\ref{eq:1.1})--(\ref{eq:1.2}) and satisfies the basic assumptions of Section \ref{s2}.


\subsection{Portfolio selection problems}

Investigations of portfolio selection problems were started in the Markowitz works \cite{Mar52,Mar59}.
These problems still play significant role in various financial decisions.
We recall that the problem is to distribute the investment capital among some ($n$) assets, i.e. one has to
define the investment shares vector $x=(x_{1}, \ldots, x_{n})^{\top}$ such that
$$
\sum \limits^{n} _{i=1}x_{i}=1,  \ x_{i} \geq 0, \ i=1, \ldots,n;
$$
the goal is to maximize the income
$$
\sum \limits^{n} _{i=1}\xi_{i} x_{i},
$$
where $\xi_{i}$ is the precise return of the $i$-th asset, whose value is supposed to be random.
One can thus take the mean variance and expected return values
$$
V(x)=\sum \limits^{n} _{i=1}\sum \limits^{n} _{j=1}c_{ij}x_{i}x_{j} \ \mbox{and} \
E(x)=\sum \limits^{n} _{i=1} m_{i}x_{i},
$$
where $c_{ij}$ and $ m_{i}$ are the corresponding covariance and mean for these random variables.
In such a way this problem involves in fact two objectives since one should minimize the mean variance (risk)
and maximize the expected return within a feasible investment share allocation; see also \cite{SQH05}.
The classical scalar portfolio selection optimization problem
consists in adding the inequality
$$
\sum \limits^{n} _{i=1} m_{i}x_{i} \geq w,
$$
where $w$ is the desired level of the expected return and
in minimizing the mean variance over all the constraints.
Some other formulations can be found e.g. in \cite{SQH05}. Note that all
the coefficients of this problem are usually inexact and non-stationary.
By replacing  the above inequality with the corresponding penalty term in the goal function
we can obtain another scalar optimization problem:
$$
\min \limits _{x} \to   \sum \limits^{n} _{i=1}\sum \limits^{n} _{j=1}c_{ij}x_{i}x_{j}+
(\tau/p) \left( w-\sum \limits^{n} _{i=1} m_{i}x_{i}  \right)_{+}^{p},
$$
subject to
$$
\sum \limits^{n} _{i=1}x_{i}=1,  \ x_{i} \geq 0, \ i=1, \ldots,n;
$$
with $\tau>0$, $p$ is  either 1 or 2. Clearly, it
 falls into into format (\ref{eq:1.1})--(\ref{eq:1.2})
and satisfies the basic assumptions of Section \ref{s2}.


\subsection{Market equilibrium models}

Let us consider a simple two-sided equilibrium market model of a
homogeneous commodity, which follows those in
\cite{Kon06,Kon07,Kon15e}.

The model involves $m$ traders and $l$ buyers of this commodity.
Each $i$-th trader has a price function $g_{i}(x_{i})$
and chooses his/her offer volume $x_{i}$ in the capacity segment $[0, \alpha _{i}]$.
Similarly, each $j$-th buyer has  a price function $h_{j}(y_{j})$
and chooses his/her bid volume $y_{j}$ in the capacity segment $[0, \beta _{j}]$.
All the price functions are supposed to be continuous.
Let $b$ denote the value of the external excess demand.
Then we can define the feasible set of offer/bid volumes
\begin{equation} \label{eq:3.3}
U=\left\{ u=(x,y) \in \mathbb{R}^{m+l}\ \vrule \
\begin{array}{c}
\sum \limits^{m} _{i=1} x_{i} - \sum \limits^{l} _{j=1} y_{j} = b, \\
x_{i} \in [0, \alpha _{i}], i=1, \dots, m; \ y_{j} \in [0, \beta_{j}], j=1, \dots, l.
\end{array}
\right\}
\end{equation}
We say that a pair $(\bar x, \bar y)$ constitutes an
{\em equilibrium point} if $(\bar x, \bar y) \in U$
and there exists a number $\bar \lambda$ such that
\begin{equation} \label{eq:3.4}
\begin{array}{l}
g_{i}(\bar x_{i}) \left\{
\begin{array}{ll}
\geq \bar \lambda & \mbox{if} \ \bar x_{i} = 0,  \\
=\bar \lambda   & \mbox{if} \ \bar x_{i} \in (0, \alpha _{i}), \\
\leq \bar \lambda  & \mbox{if} \ \bar x_{i}=\alpha _{i},
\end{array}
\right. \\
\mbox{for} \ i=1, \dots, m;
\end{array}
\quad
\begin{array}{l}
h_{j}(\bar y_{j}) \left\{
\begin{array}{ll}
\leq \bar \lambda \quad & \mbox{if} \ \bar y_{j} = 0, \\
=\bar \lambda  \quad & \mbox{if} \ \bar y_{j} \in (0, \beta_{j}), \\
\geq \bar \lambda \quad & \mbox{if} \ \bar y_{j}=\beta_{j},
\end{array}
 \right. \\
\mbox{for} \ j=1, \dots, l.
\end{array}
\end{equation}
Obviously, the number $\bar \lambda$ is the market clearing price.
In fact, the minimal offer (bid)
volumes correspond to traders (buyers) whose prices are greater
(less) than $\bar \lambda$, and the maximal offer (bid) volumes correspond
to traders (buyers) whose prices are less (greater) than $\bar \lambda$.
The prices of other participants are equal to $\bar \lambda$ and their
volumes may be arbitrary within their capacity bounds, but should be
subordinated to the balance equation. In case $l=0$ (respectively, $m=0$),
we have a market of traders (buyers) competing for shares of  the indicated
bid (offer) amount $|b|$.

It was shown in \cite{Kon06} (see also \cite{Kon07}), that each
equilibrium point $(\bar x, \bar y)$ is a solution of VI: Find $(\bar x, \bar y) \in U$ such
that
\begin{equation}\label{eq:3.5}
\sum \limits^{m} _{i=1} g_{i} (\bar x_{i}) (x_{i} - \bar x_{i}) -
\sum \limits^{l} _{j=1} h_{j} (\bar y_{j}) (y_{j} - \bar y_{j})
\geq 0 \quad \forall (x, y) \in U;
\end{equation}
and conversely, if a pair $(\bar x, \bar y)$ solves VI (\ref{eq:3.5}), (\ref{eq:3.3}), then
there exists $\bar \lambda$ such that $(\bar x, \bar y, \bar \lambda)$
satisfies (\ref{eq:3.4}).
Moreover,  we can define the function
$$
\varphi(u)=\varphi(x, y)=\sum \limits^{m} _{i=1} \mu _{i} (x_{i}) -
\sum \limits^{l} _{j=1} \eta_{j} (y_{j}),
$$
where
$$
 \mu _{i} (x_{i}) = \int \limits_{0}^{x_{i}}g_{i}(\tau)d\tau,  \ i=1, \dots, m; \ \mbox{and} \
\eta_{j} (y_{j}) = \int \limits_{0}^{y_{j}}h_{j}(\tau)d\tau,  \ j=1, \dots, l.
$$
Then, VI (\ref{eq:3.5}) is rewritten as follows:
$$
\langle \varphi'(\bar u), u - \bar u\rangle
\geq 0 \quad \forall u \in U
$$
and it yields the optimality condition for the optimization problem:
$$
\min \limits _{u \in U} \to  \varphi(u);
$$
cf. (\ref{eq:2.3}) and (\ref{eq:1.1}). By setting $n=m+l$, $x_{m+j}=-y_{j}$ for $j=1, \dots, l$
and proper modifying the bounds as indicated in Section \ref{s2}, we obtain
a particular case of problems (\ref{eq:1.1})--(\ref{eq:1.2}) and (\ref{eq:2.3}), (\ref{eq:1.1}).
The basic assumptions of Section \ref{s2} are also satisfied.


\section{Method and its convergence}\label{s4}

We now describe a two-level method of selective bi-coordinate variations (BCV for short)
for optimization problem  (\ref{eq:1.1})--(\ref{eq:1.2})
and the related  VI (\ref{eq:2.3}), (\ref{eq:1.2})
under assumptions {\bf (A1)}--{\bf (A3)}. For brevity, set
$$
g_{il}(x)=\frac{\partial f_{l}(x)}{\partial x_{i}} \ \mbox{and} \  h_{il}(x)=g_{il}(x)/a_{il},
\ \mbox{for} \ i \in I, \ l=1,2, \ldots;
$$
$\mathbb{Z}_{+}$ denotes the set of non-negative integers, and
$\pi_{V}(u)$ denotes the projection of a point $u$ on a set $V$.
Also, given a sequence $\{\varepsilon _{l}\}$ and a point $x$, let
$$
I_{l}^{-}(x)=\{ i \in I \ | \
     x_{i} \geq \alpha '_{il}+\varepsilon _{l}/a_{il} \},  \ I_{l}^{+}(x)=\{ i \in I \ | \
     x_{i} \leq \alpha ''_{il}-\varepsilon _{l}/a_{il} \}; \ l=1,2, \ldots
$$

\medskip \noindent
 {\bf Method (BCV).} \\ {\em Initialization:} Choose a point $z^{0} \in D_{0}$, numbers $\sigma \in (0,1)$,
$\theta  \in (0,1)$, and
sequences $\{\delta _{l}\} \searrow 0$, $\{\varepsilon _{l}\}
\searrow 0$. Set $l=1$.\\
{\em Step 0:} Set  $k=0$, $x^{0}=\pi_{D_{l}}(z^{l-1})$.\\
{\em Step 1:}  Choose a pair of indices $i \in I_{l}^{-}(x^{k})$ and $ j \in
I_{l}^{+}(x^{k})$ such that
\begin{equation} \label{eq:4.1}
h_{il}(x^{k})-h_{jl}(x^{k}) \geq \delta _{l},
\end{equation}
set $\gamma_{k}=\min\{a_{il}(x^{k}_{i}-\alpha '_{il}), a_{jl}(\alpha ''_{jl}-x^{k}_{j})\}$,
$i_{k}=i$, $j_{k}=j$ and go to Step 2. Otherwise
(i.e. if (\ref{eq:4.1}) does not hold for
all $i \in I_{l}^{-}(x^{k})$ and $ j \in
I_{l}^{+}(x^{k})$) set $z^{l}=x^{k}$, $l=l+1$ and go to Step 0. {\em (Restart)} \\
{\em Step 2:}  Set
$$
d^{k}_{s}= \left\{ {
\begin{array}{rl}
\displaystyle
-1/a_{sl} \quad & \mbox{if} \ s=i, \\
1/a_{sl} \quad & \mbox{if} \ s=j, \\
0 \quad & \mbox{otherwise};
\end{array}
} \right.
$$
determine $m$ as the smallest number in $\mathbb{Z}_{+}$ such that
\begin{equation} \label{eq:4.2}
 f_{l} (x^{k}+\theta ^{m}\gamma_{k} d^{k})
 \leq f _{l}(x^{k})+\sigma \theta ^{m}\gamma_{k}
  \langle f_{l}'(x^{k}),d^{k} \rangle,
\end{equation}
set $\lambda_{k}=\theta ^{m}\gamma_{k}$, $x^{k+1}=x^{k}+\lambda_{k}d^{k}$, $k=k+1$ and go to Step 1.\\
\medskip

Thus, the method has a two-level structure where each outer iteration (stage) $l$
contains some number of inner iterations in $k$
with the fixed tolerances $\delta _{l}$ and  $\varepsilon _{l}$. Completing each stage,
which is marked as restart, leads to the new approximation problem (\ref{eq:2.8})--(\ref{eq:2.9})
with decreasing of the tolerances.

Note that $i_{k}\neq j_{k} $ due to (\ref{eq:4.1}), besides, $\gamma_{k} \geq \varepsilon_{l}$
and the point $x^{k}+\gamma_{k} d^{k}$ is always feasible. Moreover, by definition,
\begin{equation} \label{eq:4.3}
\mu _{kl} =\langle f_{l}'(x^{k}),d^{k} \rangle
  =h_{j_{k},l}(x^{k})-h_{i_{k},l}(x^{k}) \leq
 -\delta_{l}<0,
\end{equation}
in (\ref{eq:4.2}). It follows that
 \begin{equation} \label{eq:4.4}
 f _{l}(x^{k+1}) \leq f_{l} (x^{k})+\sigma \lambda_{k} \mu _{kl} \leq f _{l}(x^{k})-\sigma \lambda_{k}\delta_{l}.
\end{equation}
We first justify the linesearch.


\begin{lemma} \label{lm:4.1} Suppose assumptions {\bf (A2)}--{\bf (A3)} are fulfilled.
Then the linesearch procedure in Step 2 is always finite.
\end{lemma}
{\bf Proof.}
If we suppose that the linesearch procedure is infinite, then (\ref{eq:4.2}) does not hold and
$$
(\theta ^{m}\gamma_{k})^{-1 }(f_{l} (x^{k}+\theta ^{m}\gamma_{k}
d^{k}) - f _{l}(x^{k}))>\sigma \mu _{kl},
$$
for $m \to \infty$. Hence, by taking the limit we have $    \mu _{kl}
\geq \sigma \mu _{kl}$, hence $\mu _{kl} \geq 0$, a
contradiction with $\mu _{kl}  \leq -\delta_{l}<0$. \QED

We show that each stage is well defined.


\begin{proposition} \label{pro:4.1}
Suppose assumptions {\bf (A2)}--{\bf (A3)} are fulfilled.
Then the number of iterations at each stage $l$ is finite.
\end{proposition}
{\bf Proof.}
Fix any $l$. Since the sequence $\{x^{k}\}$ is contained in the bounded
 set $D_{l}$, it has limit points.
Besides, by (\ref{eq:4.4}), we have
$$
f_{l}^{*}=\min\limits_{x \in D_{l}} f_{l}(x)\leq f_{l}(x^{k})
$$
and $f_{l}(x^{k+1})\leq f_{l}(x^{k})-\sigma \delta_{l}
\lambda_{k}$, hence
$$
\lim \limits_{k\rightarrow \infty }\lambda_{k}=0.
$$
Suppose that the sequence $\{x^{k}\}$ is infinite. Since the set $I$
is finite, there is a pair of indices $(i_{k},j_{k})=(i,j)$,
 which is repeated infinitely. Take the corresponding
subsequence $\{k_{s}\}$, then $d^{k_{s}}=\bar d$, where
$$
\bar d_{t}= \left\{ {
\begin{array}{rl}
\displaystyle
-1/a_{tl} \quad & \mbox{if} \ t=i, \\
1/a_{tl} \quad & \mbox{if} \ t=j, \\
0 \quad & \mbox{otherwise}.
\end{array}
} \right.
$$
Without loss of generality, we can
suppose that the subsequence $\{x^{k_{s}}\}$ converges to a point
$\bar x$ and due to (\ref{eq:4.3}) we have
$$
\langle f_{l}'(\bar x),\bar d \rangle= \lim \limits_{s\rightarrow \infty
} \langle f_{l}'(x^{k_{s}}),\bar d \rangle \leq -\delta_{l}.
$$
However,  (\ref{eq:4.2}) does not hold for the stepsize $\lambda_{k}/\theta$.
Setting $k=k_{s}$ gives
$$
(\lambda_{k_{s}}/\theta)^{-1 }(f_{l}
(x^{k_{s}}+(\lambda_{k_{s}}/\theta) \bar d) - f_{l}
(x^{k_{s}}))>\sigma \langle f_{l}'(x^{k_{s}}),\bar d \rangle,
$$
hence, by taking the limit $s\rightarrow \infty$ we obtain
$$
\langle f_{l}'(\bar x),\bar d \rangle= \lim \limits_{s\rightarrow \infty
} \left\{(\lambda_{k_{s}}/\theta)^{-1 }(f_{l}
(x^{k_{s}}+(\lambda_{k_{s}}/\theta) \bar d) - f_{l}
(x^{k_{s}}))\right\} \geq  \sigma \langle f'(\bar x),\bar d
\rangle,
$$
i.e.,  $   (1-\sigma ) \langle f_{l}'(\bar x),\bar d \rangle \geq 0$,
which is a contradiction. \QED

We are ready to prove convergence of the whole method.


\begin{theorem} \label{thm:4.1}
Under assumptions {\bf (A1)}--{\bf (A3)} it holds that:

(i)  the number of changes of index $k$ at each stage $l$  is finite;

(ii) the sequence $\{z^{l}\}$ generated by method (BCV) has limit points, all
these limit points are solutions of VI (\ref{eq:2.3}),
(\ref{eq:1.2});

(iii) if $f$ is semi-convex, then
\begin{equation} \label{eq:4.5}
 \lim \limits_{l\rightarrow \infty} f(z^{l})=f^{*};
\end{equation}
and all the limit points of $\{z^{l}\}$ belong to $D^{*}$.
\end{theorem}
{\bf Proof.} Assertion (i) has been obtained in Proposition \ref{pro:4.1}.
Due to assumptions {\bf (A1)}--{\bf (A2)}, the sets $\{D_{l}\}$ are uniformly bounded.
Then the sequence $\{z^{l}\}$ is bounded, hence it has limit points.
Take an arbitrary limit point $\bar z$ of $\{z^{l}\}$, then $\bar z\in D$ due to {\bf (A2)} and
$$
\lim \limits_{s\rightarrow \infty }z^{l_{s}}=\bar z.
$$
Let $p$ and $q$ be arbitrary indices such that $\bar z_{p}\in (\alpha'_{p},\alpha ''_{p}]$
and $\bar z_{q} \in [\alpha'_{q},\alpha ''_{q})$.
Then
$z^{l_{s}}_{p} \geq \alpha '_{p,l_{s}}+\varepsilon_{l_{s}}/a_{p,l_{s}}$ and
$z^{l_{s}}_{q} \leq \alpha ''_{q,l_{s}}-\varepsilon_{l_{s}}/a_{q,l_{s}}$, i.e.,
$p \in I_{l_{s}}^{-}(z^{l_{s}})$ and $q \in I_{l_{s}}^{+}(z^{l_{s}})$,
for $s$ large enough, hence
$$
  h_{p,l_{s}}(z^{l_{s}})-h_{q,l_{s}}(z^{l_{s}}) \leq \delta_{l_{s}}
$$
due to the stopping rule in Step 1.
Taking here the limit $s\rightarrow \infty$ and applying
{\bf (A2)}--{\bf (A3)}, we obtain
$$
(1/a_{p})\bar g_{p}\leq (1/a_{q})\bar g_{q}
$$
for some $\bar g \in \partial^{\uparrow}f(\bar z)$.
This means that the point $\bar z$ satisfies the optimality
conditions (\ref{eq:2.7}). Due to Proposition \ref{pro:2.2}, $\bar
z$ solves VI (\ref{eq:2.3}), (\ref{eq:1.2}) and assertion (ii) holds.

Next, if $f$ is
semi-convex, then by Proposition \ref{pro:2.1} each limit point
$\bar z=\lim \limits_{s\rightarrow \infty }z^{l_{s}}$
of $\{z^{l}\}$ solves problem (\ref{eq:1.1})--(\ref{eq:1.2}) and
$$
\lim \limits_{s\rightarrow \infty }f(z^{l_{s}})=f(\bar z)=f^{*}
$$
due to the continuity of $f$.
Since the subsequence $\{z^{l_{s}}\}$ was taken arbitrarily,
this gives (\ref{eq:4.5}) and assertion (iii).
\QED


\section{Modifications and applications}\label{s5}

First of all we would like to emphasize the fact that
convergence of the method (BCV) is attained without any concordance rules
of approximation accuracy for the problem data, cost function,
and threshold tolerances. We do not impose any condition for
approximation of solution accuracy for intermediary
problem (\ref{eq:1.1})--(\ref{eq:1.2}) or
the related  VI (\ref{eq:2.3}), (\ref{eq:1.2}). Note that any
explicit indication of this accuracy is not easy since we
do not require (strong) convexity of the cost function.

The method described admits various modifications.
We briefly discuss some of them now. Concerning the implementation of the method, we
note that utilization of the projection onto the current
feasible set $D_{l}$ in Step 0 is not obligatory.
The main condition is $x^{0} \in D_{l}$, but the other additional condition
$ f_{l} (x^{0}) \leq f _{l}(z^{l-1})$ may give better performance.

We described the method with the current type
Armijo linesearch procedure for more generality.
However, some other stepsize rules can be applied in the method
with maintaining all the results of Section \ref{s4}.
For instance, a linesearch procedure based on calculation of
only two gradient components was proposed
and substantiated in \cite{Kon16a}
for the case of convex function. Similarly, we can replace
(\ref{eq:4.2}) with the
following rule:
$$
\langle f_{l}'(x^{k}+\theta ^{m}\gamma_{k} d^{k}),d^{k} \rangle
 \leq \sigma \theta ^{m}\gamma_{k} \langle f_{l}'(x^{k}),d^{k} \rangle,
$$
or equivalently,
$$
h_{j_{k},l}(x^{k}+\theta ^{m}\gamma_{k} d^{k})-h_{i_{k},l}(x^{k}+\theta ^{m}\gamma_{k} d^{k})
 \leq \sigma \theta ^{m}\gamma_{k} (h_{j_{k},l}(x^{k})-h_{i_{k},l}(x^{k}) ).
$$
Its preference stems from the fact that the vector $d^{k}$ has only two non-zero coordinates.

We can even drop the linesearch and calculate the stepsize $\lambda_{k}$
explicitly if partial gradients of the goal function are Lipschitz continuous.
For the bi-coordinate methods these stepsize rules were
substantiated in \cite{Bec14,Kon16a}. Application of this rule
to (BCV) and substantiation can be made similarly, hence we leave
this part for the interested reader and refer to
\cite{Bec14,Kon16a} for more discussion. We only observe that
the explicit stepsize rule reduces computational expenses essentially,
but requires rather precise estimates of
the corresponding Lipschitz constants that may create difficulties in the case of
a general nonlinear cost function.

We now turn to application of the method to the market equilibrium models
from Section \ref{s2}. It was noticed in \cite{Kon16a} that the selective bi-coordinate
method proposed there can be treated as a decentralized dynamic exchange process for attaining
equilibrium states in one-sided and two-sided
markets. Each iteration is treated as a bilateral transaction
for a pair of participants (economic agents) after comparison of
their price difference in (\ref{eq:4.1}). Then the agents
simultaneously change their current transaction amounts
in order to keep the balance and bound constraints.
The agents reduce the transaction
thresholds ($\delta _{l}$ and $\varepsilon _{l}$)
sequentially if the current values appear too big (restart); see
\cite{Kon15e} for more details and comparisons.

The results of Section \ref{s3} enlarge the field of applications
of such processes essentially. In fact, it was shown in
Section \ref{s2} that the corresponding market equilibrium model involves
those in \cite{Kon16a,Kon15e} as particular cases. More precisely,
both one-sided and two-sided models from \cite{Kon16a,Kon15e}
can be written in the compact format (\ref{eq:1.1})--(\ref{eq:1.2})
or (\ref{eq:2.3}), (\ref{eq:1.2}), which gives a simpler
process definition in comparison with that in \cite[Section 6]{Kon16a}.
Besides, our current formulation now handles both upper and lower
bounds for variables.

Moreover, we note that after transformation of the market equilibrium model from
Section \ref{s2} into format (\ref{eq:2.3}), (\ref{eq:1.2})
we can differ agents by considering signs of their volume variables.
That is, $x_{i}>0$ indicates the $i$-th offer value, whereas
$x_{j}<0$ indicates the $j$-th bid value $|x_{j}|$. In the models from
\cite{Kon16a,Kon15e}, both the upper and lower
bounds of one agent must have the same sign, hence  his/her role is fixed
as either trader or buyer. However, we can now  utilize different signs for
upper and lower bounds of one agent in format (\ref{eq:2.3}), (\ref{eq:1.2}),
say, $\alpha' _{i} <0$ and $\alpha'' _{i}>0$. This means that
the $i$-th agent can change his/her role in this market model.
Therefore, the results of Section \ref{s3} confirm that
the selective bi-coordinate method proposed can serve
 as a decentralized dynamic exchange process
in much more complex and non-stationary market systems.


\section{Computational experiments}\label{s6}

In order to check the performance of the proposed method we carried
out series of computational experiments on test problems.
For comparison, we took  the known conditional gradient method
(CGM) \cite{FW56,LP66} and  marginal-based bi-coordinate descent method (MBC)
\cite{Kor80,LPR09} with the same  Armijo linesearch procedures.
We recall that the computation of the descent direction in (MBC), unlike (BCV),
is based on finding the so-called most violated pair of indices.
All the methods were implemented in {\em Delphi} with double precision
arithmetic. The main goal was to compare convergence of the methods
despite the smaller iteration expenses of (BCV). In all the cases,
we took the accuracy $\mu=0.1$ and the starting point $(\beta/n)e$,
where $e$ denote the vector of units in $\mathbb{R}^{n} $. We chose
$\sigma =\theta =0.5$, and the rule $\delta_{l+1}=\nu \delta_{l}$,
$\varepsilon_{l+1}=\nu\varepsilon_{l}$  with $\nu = 0.5$ for (BCV).
For testing, we chose problems of form (\ref{eq:1.1})--(\ref{eq:1.2}) with $a_{i}=1$,
$\alpha' _{i}=0$ and $\alpha'' _{i}=1+ (\beta/n)+ 0.5 \sin (i)$
for $i =1, \ldots, n$.

In the first two series, we took the stationary problem (\ref{eq:1.1})--(\ref{eq:1.2})
with the fixed data and smooth goal function $f$. We hence took the value
$$
\Delta (x)=\max _{y\in D}\langle f'(x),x-y \rangle
$$
as an error bound at $x$ and write $\Delta_{k}$ for the accuracy $\Delta(x)$
after full $k$ iterations.
We give the total number of iterations of each method for attaining the indicated accuracy
in each case, sign \lq\lq -" means that the error was too big, namely, $\Delta_{500}>1$.

In the first series, we took the quadratic cost function $f(x)=\varphi(x)$ where
$$
\varphi (x)= 0.5 \langle Px,x \rangle,
$$
the elements of the matrix $P$ were defined by
$$
p_{ij}= \left\{ {
\begin{array}{rl}
\displaystyle
\sin (i) \cos (j) \quad & \mbox{if} \ i<j, \\
\sin (j) \cos (i) \quad & \mbox{if} \ i>j, \\
\sum \limits_{i=1}^{n} | p_{ij}| +1 \quad & \mbox{if} \ i=j.
\end{array}
} \right.
$$
We varied the parameter $\beta$ and dimensionality $n$.
The results are given in Table \ref{tbl:1}.
\begin{table}
\caption{} \label{tbl:1}
\begin{center}
\begin{tabular}{|l|c|c|c|}
\hline
                    &  (CGM)  & (BCV) & (MBC)    \\
\hline
   $\beta=5$            &    &   &     \\
\hline
    $n=10$            & 66  &  30 &  $\Delta_{500}\approx 1.28$ \\
    $n=20$            & 22  &  41 &  $\Delta_{500}\approx 0.99$ \\
     $n=50$           & 82  &  96 &  $\Delta_{500}\approx 0.91$ \\
     $n=100$          &  $\Delta_{500}\approx 0.1$     & 213 &  $\Delta_{500}\approx 1.63$   \\
\hline
   $\beta=10$            &    &   &     \\
\hline
    $n=10$            & 55  &  40 &  $\Delta_{500}\approx 5.13$ \\
    $n=20$            & 103  &  54 &  - \\
     $n=50$           & 90  &  145 &  - \\
     $n=100$          &  $\Delta_{500}\approx 0.48$     & 299 &  -   \\
\hline
   $\beta=20$            &    &   &     \\
\hline
    $n=10$            & $\Delta_{500}\approx 0.14$   &  62 &  - \\
    $n=20$            & $\Delta_{500}\approx 0.23$   &  80 &  - \\
     $n=50$           & $\Delta_{500}\approx 0.21$   &  191 &  - \\
     $n=100$          &  $\Delta_{500}\approx 1.07$     & 405 &  -   \\
\hline
\end{tabular}
\end{center}
\end{table}

In the second series, we took the convex cost function
$$
f(x)=\varphi(x)+\psi(x),
$$
where the function $\varphi$ was defined as above,
$$
\psi(x)=-\ln(\langle c,x\rangle+\xi),
$$
the elements of the vector $c$ are defined by
$$
c_{i}= 2+\sin(i) \ \mbox{ for } \ i=1,\ldots, n,
$$
and $\xi=5$.  The results are
given in Table \ref{tbl:2}.

\begin{table}
\caption{} \label{tbl:2}
\begin{center}
\begin{tabular}{|l|c|c|c|}
\hline
                    &  (CGM)  & (BCV) & (MBC)    \\
\hline
   $\beta=5$            &    &   &     \\
\hline
    $n=10$            & 77  &  29 &  $\Delta_{500}\approx 1.29$ \\
    $n=20$            & 30  &  35 &  - \\
     $n=50$           & 111  & 109 &  - \\
     $n=100$          & 457  & 240 &  -   \\
\hline
   $\beta=10$            &    &   &     \\
\hline
    $n=10$            & 62  &  44 &  $\Delta_{500}\approx 5.14$ \\
    $n=20$            & 77  &  53 &  - \\
     $n=50$           & 115  &  167 &  - \\
     $n=100$          &  $\Delta_{500}\approx 0.46$     & 282 &  -   \\
\hline
   $\beta=20$            &    &   &     \\
\hline
    $n=10$            & $\Delta_{500}\approx 0.12$   &  68 &  - \\
    $n=20$            & $\Delta_{500}\approx 0.21$   &  75 &  - \\
     $n=50$           & $\Delta_{500}\approx 0.24$   &  220 &  - \\
     $n=100$          &  $\Delta_{500}\approx 1.07$     & 350 &  -   \\
\hline
\end{tabular}
\end{center}
\end{table}

In the third series, we took the non-smooth convex cost function
$$
f(x)=\varphi(x)+\psi(x)+\sum \limits^{n} _{i=1} |x_{i}|,
$$
where the functions $\varphi$ and $\psi$  were defined as above.
We also took the fixed coefficients $a_{i}=1$,
$\alpha' _{i}=0$ and $\alpha'' _{i}=1+ (\beta/n)+ 0.5 \sin (i)$
for $i =1, \ldots, n$. We utilized the smooth approximations of the form
$$
\phi(x,\tau)=\varphi(x)+\psi(x)+\sum \limits^{n} _{i=1} \sqrt{x_{i}^{2}+\tau^{2}}.
$$
In other words, we replace (\ref{eq:1.1})--(\ref{eq:1.2}) with the sequence of the following
smooth optimization problems
$$
 \min \limits _{x \in D} \to f_{l}(x),
$$
where $f_{l}(x)=\phi(x,\tau_{l})$ for some sequence $\{\tau _{l}\} \searrow 0$, i.e. set $D_{l}=D$.
The main goal was to check the performance for such
smooth approximations of the non-smooth initial problem.
Since (MBC) appeared rather slow, we compared only (CGM) and (BCV).
We used the value
$$
\Delta (x,\tau)=\max _{y\in D}\langle \phi'(x,\tau),x-y \rangle
$$
as an error bound at $x$. We stopped the calculations under the condition
$$
\Delta (x,\tau) \leq \mu \ \mbox{and} \ \tau \leq \mu
$$
with $\mu = 0.1$. We used the rule $\tau_{l+1}=\max\{\mu, \nu \tau_{l}\}$ with $\nu = 0.5$,
the other parameters of the methods were chosen as above.
The results are given in Table \ref{tbl:3}, where $\Delta_{\tau,k}$
denotes the accuracy $\Delta(x,\tau)$
after full $k$ iterations.

\begin{table}
\caption{} \label{tbl:3}
\begin{center}
\begin{tabular}{|l|c|c|}
\hline
                    &  (CGM)  & (BCV)    \\
\hline
   $\beta=5$            &    &     \\
\hline
    $n=10$            & 154  &  57  \\
    $n=20$            & 43 &  52  \\
     $n=50$           & 103  & 85  \\
     $n=100$          & 300  & 234   \\
\hline
   $\beta=10$            &    &   \\
\hline
    $n=10$            & 98  &  49  \\
    $n=20$            & 123  &  52 \\
     $n=50$           & 183  &  136 \\
     $n=100$          &  $\Delta_{\tau,500}\approx 0.81, \tau=0.8$     & 271   \\
\hline
   $\beta=20$            &    &       \\
\hline
    $n=10$            & $\Delta_{\tau,500}\approx 6.8, \tau=6.4$  &  66  \\
    $n=20$            & $\Delta_{\tau,500}\approx 0.9, \tau=0.8$   &  67  \\
     $n=50$           & $\Delta_{\tau,500}\approx 0.48, \tau=0.2$   &  197  \\
     $n=100$          &  $\Delta_{\tau,500}\approx 2.06, \tau=1.6$    & 468    \\
\hline
\end{tabular}
\end{center}
\end{table}
In all the experiments, (BCV) showed rather rapid and stable convergence.
In almost all the cases, (BCV) showed better results than (CGM).
At the same time, these experiments showed rather slow and instable convergence of
(MBC). We also noticed that the presence of nonlinear  functions
or approximations of non-smooth functions had no significant influence on the
convergence of (BCV).


\section{Conclusions}

We suggested a new method of bi-coordinate variations for
non-stationary and non-smooth optimization problems,
which involve two side constraints for variables  and a single linear equality.
Its descent direction rule selects
only two coordinates for changes and enables us to avoid
calculation of all the gradient components of a current smooth approximation function
at each iteration. Therefore, the new method is simpler essentially than the usual
gradient or dual type ones, but  does not impose any concordance rules
of approximation accuracy for the problem data, cost function,
and threshold tolerances. We showed some its possible fields of applications.
Its convergence was established under rather mild assumptions.
Computational tests showed certain preferences of the proposed method over the
known ones.

\section*{Acknowledgement}

This work was supported by the RFBR grant, project No. 16-01-00109a
and by grant No. 297689 from Academy of Finland.


\end{document}